\documentclass[12pt,leqno]{article}
\usepackage[doc]{optional}
\usepackage{color}
\usepackage{float}
\usepackage{soul}
\usepackage{graphicx}

\definecolor{labelkey}{rgb}{0,0.08,0.45}
\definecolor{refkey}{rgb}{0,0.6,0.0}

\definecolor{Brown}{rgb}{0.45,0.0,0.05}
\definecolor{lime}{rgb}{0.00,0.8,0.0}
\definecolor{lblue}{rgb}{0.5,0.5,0.99}

\usepackage{multirow}
\usepackage{mathpazo}
\usepackage{enumitem}
\usepackage{empheq}

\usepackage{fancyhdr}
\fancyheadoffset[R]{0.5cm}
\usepackage{amsmath}

\usepackage{amssymb}
\usepackage{theorem}

\usepackage[margin=1in]{geometry}

\definecolor{myblue}{rgb}{.9, .9, 1}

\newcommand{\sepp}{\setlength{\itemsep}{-2pt}}

\newcommand{\nnn}{\ensuremath{{n\in{\mathbb N}}}}

\newcommand{\menge}[2]{\big\{{#1}~\big |~{#2}\big\}}

\newcommand{\fenv}[1]%
{\ensuremath{\,\overrightarrow{\operatorname{env}}_{#1}}}
\newcommand{\benv}[1]%
{\ensuremath{\,\overleftarrow{\operatorname{env}}_{#1}}}
\newcommand{\emp}{\ensuremath{\varnothing}}
\newcommand{\ve}{\ensuremath{\varepsilon}}

\newcommand{\scal}[2]{\left\langle{#1},{#2}  \right\rangle}

\newcommand{\RR}{\ensuremath{\mathbb R}}
\newcommand{\RP}{\ensuremath{\mathbb{R}_+}}
\newcommand{\RPP}{\ensuremath{\mathbb{R}_{++}}}
\newcommand{\RM}{\ensuremath{\mathbb{R}_-}}

\newcommand{\RX}{\ensuremath{\,\left]-\infty,+\infty\right]}}

\newcommand{\reli}{\ensuremath{\operatorname{ri}}}
\newcommand{\inte}{\ensuremath{\operatorname{int}}}

\newcommand{\closu}{\ensuremath{\operatorname{cl}}}

\newcommand{\card}{\ensuremath{\operatorname{card}}}

\newcommand{\conv}{\ensuremath{\operatorname{conv}}}

\newcommand{\aff}{\ensuremath{\operatorname{aff}}}

\newcommand{\pinf}{\ensuremath{+\infty}}

{\begin{list}{}{%
\settowidth{\labelwidth}{\textrm{#1~}}%
\setlength{\leftmargin}{\labelwidth+\labelsep}}}
{\end{list}}
\newtheorem{theorem}{Theorem}[section]
\newtheorem{lemma}[theorem]{Lemma}
\newtheorem{corollary}[theorem]{Corollary}
\newtheorem{proposition}[theorem]{Proposition}
\newtheorem{definition}[theorem]{Definition}
\theoremstyle{plain}{\theorembodyfont{\rmfamily}
}
\theoremstyle{plain}{\theorembodyfont{\rmfamily}
}
\theoremstyle{plain}{\theorembodyfont{\rmfamily}
}
\theoremstyle{plain}{\theorembodyfont{\rmfamily}
\newtheorem{example}[theorem]{Example}}
\newtheorem{fact}[theorem]{Fact}
\theoremstyle{plain}{\theorembodyfont{\rmfamily}
\newtheorem{remark}[theorem]{Remark}}
\def\proof{\noindent{\it Proof}. \ignorespaces}
\def\endproof{\ensuremath{\hfill \quad \square}}

\allowdisplaybreaks 

\def\tto{\rightrightarrows}
\def\ox{\overline{x}}
\def\oy{\overline{y}}

\def\O{\Omega}
\def\mcF{{\mathcal{F}}}
\def\DM{{\rm D}}

\def\mcA{{\mathcal{A}}}

\def\disp{\displaystyle}
\begin{document}

\title{On the convexity of piecewise-defined functions}

\author{
Heinz H.\ Bauschke\thanks{ 
Mathematics, University of British Columbia, Kelowna, B.C.\ V1V~1V7, Canada. E-mail: \texttt{heinz.bauschke@ubc.ca}.},
Yves Lucet\thanks{Computer Science, University of British Columbia, Kelowna, B.C.\ V1V~1V7, Canada. E-mail: \texttt{yves.lucet@ubc.ca}},~and~
Hung M.\ Phan\thanks{ Mathematics, University of British Columbia, Kelowna, B.C.\ V1V\,1V7, Canada.
E-mail: \texttt{hung.phan@ubc.ca}.}}
\date{August 16, 2014}

\maketitle

%
\begin{abstract} \noindent
Functions that are piecewise defined are a common sight in
mathematics while convexity is a property
especially desired in optimization. 
Suppose now a piecewise-defined function
is convex on each of its defining components --- when can we conclude that
the entire function is convex? In this paper we provide several
convenient, verifiable conditions guaranteeing convexity (or the
lack thereof). Several examples are presented to illustrate
our results.
\end{abstract}

{\small

}
\section{Introduction}

Consider the function 
\begin{equation}\label{e:0715a}
f\colon \RR^2\to\RR \colon (x,y)\mapsto \begin{cases}
\disp\frac{x^2+y^2+2\max\{0,xy\}}{|x|+|y|},&\text{if $(x,y)\neq
(0,0)$;}\\
0,&\text{otherwise.}
\end{cases}
\end{equation}
Clearly, $f$ is a piecewise-defined function 
with continuous components
\begin{subequations}\label{e:0715b}
\begin{alignat}{3}
& f_1(x,y) &&:=x+y &&\qquad\text{on}\qquad A_1:=\RP\times\RP;\\
& f_2(x,y) &&:=\frac{x^2+y^2}{-x+y} &&\qquad\text{on}\qquad A_2:=\RM\times\RP;\\
& f_3(x,y) &&:=x+y&&\qquad\text{on}\qquad A_3:=\RM\times\RM;\\
& f_4(x,y) &&:=\frac{x^2+y^2}{x-y} &&\qquad\text{on}\qquad A_4:=\RP\times\RM.
\end{alignat}
\end{subequations}
One may check that each $f_i$ is a convex function 
(see Example~\ref{ex:1} below for details). 
However, \emph{whether or not $f$ itself is convex} is not
immediately clear. (As it turns out, $f$ is convex.)

On the other hand, if 
\begin{subequations}
\begin{align}
f_1(x)=x\quad\text{on}\quad A_1:= \RM;\\
\text{and}\quad
f_2(x)=-x\quad\text{on}\quad A_2:=\RP,
\end{align}
\end{subequations}
then $f_1$ and $f_2$ are convex while 
 the induced piecewise-defined function $f(x)=-|x|$
is \emph{not} convex.

These and similar examples motivate the goal of this paper which
is to present 
{\em verifiable conditions guaranteeing the convexity of a
piecewise-defined function provided that each component is
convex.} 
Special cases of our results have been known in the convex
interpolation community (see Remark~\ref{r:ci}). 
Moreover, our results have applications to computer-aided
convex analysis (see Remark~\ref{r:plq}).

The remainder of this paper is organized as follows.
In Section~\ref{s:convfunc}, we collect various auxiliary results
concerning convexity and differentiability.
We also require properties of collections of sets and of
functions which we develop in Section~\ref{s:Compsets}
and Section~\ref{s:Compfunc}, respectively. 
Our main results guaranteeing convexity are presented in
Section~\ref{s:main}. 
Various examples illustrating convexity and the lack thereof
are discussed in Section~\ref{s:Checkconv} and 
Section~\ref{s:Lackconv}, respectively.

\noindent{\bf Notation:}
Throughout, $X$ is a Euclidean space
with inner product $\scal{\cdot}{\cdot}$ and
induced norm $\|\cdot\|$.
$\RR$ denotes the set of real numbers, $\RP:=\menge{x\in\RR}{x\geq0}$,
and $\RM:=-\RP$.  For $x$ and $y$ in $X$, $[x,y]:=\menge{(1-t)x+ty}{0\leq
t\leq 1}$ is the line segment connecting $x$ and $y$. Similarly,
we set $\left]x,y\right[:=\menge{(1-t)x+ty}{0<t<1}$,
$\left[x,y\right[:=\menge{(1-t)x+ty}{0\leq t<1}$, and
$\left]x,y\right]:=\left[y,x\right[$. For a subset $A$ of $X$,
$\conv A$, $\closu A$, $\inte A$, $\aff A$, and $\reli A$ respectively
denote the \emph{convex hull}, the \emph{closure}, 
the \emph{interior}, the \emph{affine hull},
and the \emph{relative interior of $A$}. 
Furthermore, $\iota_A$ is the \emph{indicator function}
of $A$ defined by $\iota_A(x)=0$, if $x\in A$; and $+\infty$ otherwise.
Let $f\colon X\to\RX = \RR\cup\{\pinf\}$. 
The \emph{domain} of $f$ is $\DM_f:=\menge{x\in
X}{f(x)<+\infty}$; $f$ is said to be {\em proper} if $\DM_f\neq\emp$.
The restriction of $f$ on some subset $A$ of $X$ is denoted by
$f\big|_{A}$. A set-valued mapping $F$ from $X$ to another Euclidean
space $Y$ is denoted by $F:X\tto Y$; and its domain is 
$\DM_F:=\menge{x\in X}{F(x)\neq\emp}$.
For further background and notation, we refer the reader to
\cite{BC2011, MorNam, ROCKAFELLAR-70a, ROCKAFELLAR-98a, ZALINESCU-02}.

\section{Convexity and differentiability}

\label{s:convfunc}

Let $f:X\to\RX$ be proper. For every $x\in X$, the
\emph{subdifferential} (in the sense of convex analysis) of $f$ at $x$, 
denoted by $\partial
f(x)$, is the set of all vectors $x^*\in X$ such that
\begin{equation}
(\forall y\in X)\quad\scal{x^*}{y-x}\leq f(y)-f(x).
\end{equation}
The induced operator  $\partial f:X\tto X$ has domain 
$\DM_{\partial f}=\menge{x\in X}{\partial f(x)\neq\emp}\subseteq \DM_f$.

Let us now present some auxiliary results concerning the 
convexity of a function.

\begin{lemma}
\label{l:0603a}
Let $f:X\to\RX$ and let $z_1$ and $z_2$ be in $\DM_f$.
Set $x:=(1-t)z_1+tz_2$, where $t\in[0,1]$, 
and assume that $x\in \DM_{\partial f}$. 
Then $f(x)\leq (1-t)f(z_1)+tf(z_2)$.
\end{lemma}
\proof 
Let $x^*\in\partial f(x)$. 
Then $\scal{x^*}{z_1-x}\leq f(z_1)-f(x)$
\ and \
$\scal{x^*}{z_2-x}\leq f(z_2)-f(x)$.
Hence
\begin{subequations}
\begin{align}
(1-t)\scal{x^*}{z_1-x}&\leq (1-t)(f(z_1)-f(x));\\
t\scal{x^*}{z_2-x}&\leq t(f(z_2)-f(x)).
\end{align}
\end{subequations}
Adding up the last two inequalities, 
we obtain $0\leq (1-t)f(z_1)+tf(z_2)-f(x)$.
\endproof


\begin{fact}
{\rm (See \cite[Theorem~6.1--6.3]{ROCKAFELLAR-70a}.)}
\label{f:ri}
Let $A$ be a nonempty convex subset of $X$. 
Then the following hold:
\begin{enumerate}
\item\label{f:ri-i} 
$(\forall x\in\closu A)(\forall y\in \reli A)$ 
$\left]x,y\right] \subseteq \reli A$.
\item\label{f:ri-ii}
$\reli A$ is nonempty and convex.
\item\label{f:ri-iii}
$\closu(\reli A)=\closu A$.
\end{enumerate}
\end{fact}


\begin{fact}
{\rm (See \cite[Theorem~23.4]{ROCKAFELLAR-70a}.)}
\label{f:0605a}
Let $f\colon X\to\RX$ be convex and proper.
Then $\reli \DM_f \subseteq \DM_{\partial f}$.
\end{fact}


\begin{fact}{\rm (See \cite[Theorem~2.4.1(iii)]{ZALINESCU-02})}
\label{f:0417b}
Let $f\colon X\to\RX$ be proper.
Assume that $\DM_f = \DM_{\partial f}$ is convex. 
Then $f$ is convex.
\end{fact}
\proof
Take $z_1$ and $z_2$ in $\DM_f$ and let $t\in\left[0,1\right]$.
Set
$x:=(1-t)z_1+t z_2$. Since $\partial f(x)\neq\emp$, 
Lemma~\ref{l:0603a} implies that 
$f(x)\leq(1-t)f(z_1)+tf(z_2)$. 
Therefore, $f$ is convex.
\endproof

In the presence of continuity, Fact~\ref{f:0417b} admits the
following extension.

\begin{lemma}\label{l:0416a}
Let $f\colon X\to\RX$ be proper.
Assume that 
$\DM_f$ is convex, 
that $f\big|_{ \DM_f}$ is continuous, 
and that 
$\reli \DM_f \subseteq \DM_{\partial f}$. 
Then $f$ is convex.
\end{lemma}
\proof
Fact~\ref{f:ri} implies that $\reli \DM_f$ is convex and
$\closu(\reli \DM_f)=\closu \DM_f$. Then the function
$f+\iota_{\reli \DM_f}$ is convex by Fact~\ref{f:0417b}. Since
$f\big|_{\DM_f}$ is continuous and $\reli \DM_f$ is dense in
$\DM_f$, we conclude that $f\big|_{\DM_f}$ is convex.
\endproof

Given a nonempty subset $A$ of $X$,
we define the dimension of $A$ to be 
the dimension of the linear subspace parallel to
the affine hull of $A$, i.e., 
$\dim A :=  \dim(\aff A - \aff A)$. 
We then have the following result.

\begin{lemma}\label{l:0513c}
Let $f\colon X\to\RX$ be proper.
Assume that $f\big|_{\DM_f}$ is continuous,
that $\DM_f$ is convex and at least $2$-dimensional,
and that there exists a finite subset $E$ of $X$ 
such that $f\big|_{[x,y]}$ is convex for every segment $[x,y]$
contained in $(\reli \DM_f)\smallsetminus E$. 
Then $f$ is convex.
\end{lemma}
\proof
Take two distinct points $x$ and $y$ in $\DM_f$, let 
$t\in[0,1]$, and set $z:=(1-t)x+t y$. Then 
$z\in \DM_f$ because $\DM_f$ is convex.
It remains to show that 
\begin{equation}\label{e:l0513c1}
f(z)\leq (1-t)f(x)+ t f(y).
\end{equation}
First, since $\DM_f$ is convex and $\dim(\DM_f) \geq 2$, 
there exists 
$w\in(\reli \DM_f)\smallsetminus\aff\{x,y\}$.
For each $\ve\in\left]0,1\right[$,
set 
\begin{equation}
x_\ve:=(1-\ve)x+\ve w\ ,\
y_\ve:=(1-\ve)y+\ve w\ ,\ \text{and}\ \
z_\ve:=(1-t)x_\ve+t y_\ve.
\end{equation}
Using 
Fact~\ref{f:ri}\ref{f:ri-i} and the finiteness of $E$,
we have 
\begin{equation}
z_\ve \in [x_{\ve},y_{\ve}]\subseteq (\reli \DM_f)\smallsetminus
E,
\quad
\text{for every $\ve>0$ sufficiently small.}
\end{equation}
Because $f$ is convex on $[x_{\ve},y_{\ve}]$, this implies
\begin{equation}
f(z_\ve)\leq (1-t)f(x_\ve)+ t f(y_\ve).
\end{equation}
Letting $\ve\to 0^+$, we obtain \eqref{e:l0513c1} by using the continuity of $f\big|_{ \DM_f}$.
\endproof`

\begin{remark}[the assumption on the dimension is important]
Lemma~\ref{l:0513c} fails on $\RR$ in the following sense.
Consider $f\colon\RR\to\RR\colon x\mapsto-|x|$ and set $E:=\{0\}$. 
Then all assumptions of Lemma~\ref{l:0513c} hold except that $\DM_f=\RR$ 
is only $1$-dimensional. Clearly, the conclusion of Lemma~\ref{l:0513c} 
is not true because $f$ is not convex.
\end{remark}

We now turn our attention to differentiability properties. 
Recall that $f\colon X\to\RX$ is differentiable at 
$x\in\inte \DM_f$ if there exists $\nabla f(x)\in X$ such that
$(\forall y\in X)$ $f(y)-f(x)-\scal{\nabla f
(x)}{y-x}=o(\|y-x\|)$;
$f$ is differentiable on subset $A$ of $\inte\DM_f$ if $f$ 
is differentiable at every $x\in A$.
We will require the following results.

\begin{fact}
\label{f:r25-1}
{\rm(See \cite[Theorem~25.1]{ROCKAFELLAR-70a}.)}
Let $f\colon X\to\RX$ be convex and proper,
and assume that $x\in\inte\DM_f$. 
Then $f$ is differentiable at $x$ if and only if  
$\partial f(x)$ is a singleton.
\end{fact}

\begin{fact}
\label{f:r25-5}
{\rm(See \cite[Theorem~25.5]{ROCKAFELLAR-70a}.)}
Let $f\colon X\to\RX$ be convex and proper, and 
let $\O$ be the set of points where $f$ is differentiable. 
Then $\O$ is a dense subset of $\inte \DM_f$, and its complement in 
$\inte \DM_f$ is a set of measure zero. 
Moreover, $\nabla f\colon \O\to X$ is continuous.
\end{fact}

\begin{fact}
\label{f:r25-6}
{\rm(See \cite[Theorem~25.6]{ROCKAFELLAR-70a}.)}
Let $f\colon X\to\RX$ be convex and proper such that 
$\DM_f$ is closed with nonempty interior. Then
\begin{equation}
(\forall x\in \DM_f)\quad\partial f(x)=\closu(\conv S(x)) + N_{\DM f}(x),
\end{equation}
where $N_{\DM f}(x):=\menge{x^*\in X}{(\forall y\in\DM_f)\
\scal{x^*}{y-x}\leq 0}$ is the normal cone to $\DM_f$ at $x$ and
$S(x)$ is the set of all limits of sequences $(\nabla
f(x_n))_{\nnn}$
such that $f$ is differentiable at every $x_n$ and $x_n\to x$.
\end{fact}
\section{Compatible systems of sets}

\label{s:Compsets}

In this section, we always assume that 
\begin{subequations}
\label{e:AA}
\begin{align}
&\text{$I$ is a nonempty finite set;}\\
&\text{$\mcA:=\{A_i\}_{i\in I}$ \ is a system of convex subsets of $X$;}\\
&
A:=\textstyle\bigcup_{i\in I}A_i.
\end{align}
\end{subequations}

\begin{definition}[compatible systems of sets]
\label{d:comp} 
Assume \eqref{e:AA}. 
We say that $\mcA$ is a \emph{compatible system of sets} if
\begin{equation}
\left.
\begin{matrix}
i\in I\\
j\in I\\
i\neq j
\end{matrix}
\right\}
\quad
\Rightarrow
\quad
\closu A_i\cap\closu A_j\cap\reli A=A_i\cap A_j\cap\reli A;
\end{equation}
otherwise, we say that $\mcA$ is {\em incompatible}.
\end{definition}


\begin{example}
\label{ex:0815}
Every system of finitely many closed convex subsets of $X$
is compatible.
\end{example}


\begin{example}[incompatible systems]
Suppose that $X=\RR^2$,
that $I=\{1,2\}$,
that $A_1=\left]0,1\right]\times[0,1]$, and 
that $A_2=[-1,0]\times[0,1]$. 
Then $A=A_1\cup A_2=[-1,1]\times[0,1]$ and 
$\reli A=\left]-1,1\right[\times\left]0,1\right[$. 
Thus, $\mcA = \{A_1,A_2\}$ is incompatible because
\begin{equation}
\closu A_1\cap\closu A_2\cap\reli A=\{0\}\times\left]0,1\right[\neq \emp=A_1\cap A_2\cap\reli A.
\end{equation}
\end{example}


\begin{definition}[colinearly ordered tuple]
\label{d:colot}
The tuple of vectors $(x_0,\ldots,x_n)\in X^n$ is 
said to be \emph{colinearly ordered} if 
the following hold:
\begin{enumerate}
\item $[x_0,x_n]=[x_0,x_1]\cup\cdots\cup[x_{n-1},x_n]$;
\item $0\leq\|x_0-x_1\|\leq\|x_0-x_2\|\leq\cdots\leq \|x_0-x_n\|$.
\end{enumerate}
\end{definition}




\begin{proposition}
\label{p:colin}
Assume \eqref{e:AA} and that $\mcA$ is a compatible system of sets 
(recall Definition~\ref{d:comp}). 
Then for every segment $[x,y]$ contained in $\reli A$, 
there exists a colinearly ordered tuple 
$(x_0,\ldots,x_n)$ and $\{A_{i_1},\ldots,A_{i_n}\}\subseteq \mcA$ 
such that
\begin{equation}\label{e:p-coli}
x_0=x;\quad x_n=y;\quad\text{and}\quad
\big(\forall k\in\{1,\ldots,n\}\big)\quad [x_{k-1},x_k]\subseteq A_{i_k}.
\end{equation}
\end{proposition}
\proof
Let $[x,y]\subseteq \reli A$, with $x\in A_{i_1}$ for some $i_1\in I$. 
Set $x_0:=x$. For every $t\in[0,1]$, define
\begin{equation}
x(t):=(1-t)x+ty.
\end{equation}
Furthermore, set 
\begin{equation}
t_1:=\sup\menge{t\in[0,1]}{x(t)\in A_{i_1}}
\quad\text{and}\quad
x_1:=x(t_1).
\end{equation}
Then $\left[x_0,x_1\right[\subseteq A_{i_1}$ and $x_1\in\closu A_{i_1}$.
Note also that $x_1\in\reli A$.

{\em Case~1:} $t_1=1$. Then $x_1=y\in\closu A_{i_1}$. Suppose that $y\not\in A_{i_1}$. Then, $y\in A_{i_2}$ for some $i_2\in I$. 
It follows that $y\in\closu A_{i_1}\cap\closu A_{i_2}\cap\reli A
=A_{i_1}\cap A_{i_2}\cap \reli A$, which is a contradiction. 
Therefore, $y\in A_{i_1}$ and we are done because 
$[x,y]\subseteq A_{i_1}$.

{\em Case~2:} $t_1<1$. Then there exist 
$\ve\in\left]0,1-t_1\right]$ 
and $A_{i_2}\in\mcA\smallsetminus \{A_{i_1}\}$ such that
\begin{equation}
\left]x_1,x_2\right]\subseteq A_{i_2}
\quad\text{where}\quad x_2:=x(t_1+\ve).
\end{equation}
Hence $x_1\in\closu A_{i_2}$. We then have
\begin{equation}
x_1\in\closu A_{i_1}\cap\closu A_{i_2}\cap\reli A=
A_{i_1}\cap A_{i_2}\cap\reli A.
\end{equation}
So we have split $[x,y]$ into two line segments
\begin{equation}
[x_0,x_1]\subseteq A_{i_1}\cap\reli A
\quad\text{and}\quad
[x_1,y]\subseteq \big(\textstyle\bigcup_{i\in I\smallsetminus\{i_1\}}A_i\big)\cap\reli A.
\end{equation}
Next, we repeat the above process for the segment $[x_1,y]$.
Since $\mcA$ is finite, we eventually obtain \eqref{e:p-coli}.
\endproof

\begin{remark}[closedness is not necessary for compatibility]
We note that there are compatible systems of sets that are not
closed. For example, suppose that $X=\RR^2$,
that $I=\{1,2\}$,
that $A_1=\left[0,1\right[\times[0,1]$, and
that $A_2=\left]-1,0\right]\times[0,1]$. 
Then neither $A_1$ nor $A_2$ is closed. 
However, since
\begin{equation}
\closu A_1\cap \closu A_2=A_1\cap A_2=\{0\}\times[0,1],
\end{equation}
we deduce that $\mcA = \{A_1,A_2\}$ is compatible. 
\end{remark}


\begin{definition}[active index set]
Assume \eqref{e:AA}. 
For every $x\in X$, we define the 
\emph{active index set} associated with $\mcA$ by 
\begin{equation}
I_\mcA(x):=\menge{i\in I}{x\in A_i},
\end{equation}
and we will write $I(x)$ if there is no cause for confusion.
\end{definition}

\begin{proposition}
\label{p:sigle}
Assume \eqref{e:AA} and that $\mcA$ is a compatible system of sets (see
Definition~\ref{d:comp}). 
Suppose that $x\in\inte A$ and that 
$I_\mcA(x) = \{i\}$.
Then $x\in\inte A_i$.
\end{proposition}
\proof Because $\inte A\neq\emp$, we have $\reli A=\inte A$.
Suppose to the contrary that $x\not\in\inte A_i$. 
Then there exist $j\in I\smallsetminus\{i\}$ and a 
sequence $(x_n)_\nnn$ in $A_j$ such that 
that $x_n\to x$. 
It follows that
\begin{equation}
x\in\closu A_i\cap\closu A_j\cap\reli A=A_i\cap A_j\cap\reli A,
\end{equation}
which is absurd because $I(x)=\{i\}$ by assumption. 
\endproof

\section{Compatible systems of functions}

\label{s:Compfunc}

In this section, we always assume that 
\begin{subequations}
\label{e:BB}
\begin{align}
&\text{$I$ is a nonempty finite set;}\\
&\text{$\mcF:=\{f_i\}_{i\in I}$ \ is a system of proper convex
functions from $X$ to $\RX$;}\\
&
\text{$f := \textstyle \min_{i\in I} f_i$ is the piecewise-defined function
associated with $\mcF$;}\\
&
\text{$I_{\mcF}\colon X\to I \colon x \mapsto \menge{i\in I}{x\in
\DM_{f_i}}$ is the active index set function.}
\end{align}
\end{subequations}
We will write $I(x)$ instead of $I_{\mcF}(x)$ if there is no
cause for confusion. 
Note that  $\DM_f=\bigcup_{i\in I}\DM_{f_i}$.

\begin{definition}[compatible systems of functions]
\label{d:fcomp} 
Assume \eqref{e:BB}. 
We say that $\mcF$ is a \emph{compatible system of functions} if
$(\forall i\in I)$ 
$f_i\big|_{D_{f_i}}$ is continuous and
\begin{equation}
\left.
\begin{matrix}
i\in I\\
j\in I\\
i\neq j
\end{matrix}
\right\}
\quad
\Rightarrow
\quad
f_i\big|_{\DM_{f_i}\cap \DM_{f_j}}\equiv f_j\big|_{\DM_{f_i}\cap
\DM_{f_j}}.
\end{equation}
\end{definition}

We start with a useful lemma.


\begin{lemma}\label{l:0502a}
Assume \eqref{e:BB} and that $\mcF$ 
is compatible system of functions (recall
Definition~\ref{d:fcomp}). 
Then
\begin{equation}\label{e:l0502a1}
(\forall x\in X)\quad\partial f(x)\subseteq\bigcap_{i\in I_\mcF(x)}\partial f_i(x).
\end{equation}
\end{lemma}
\proof 
Suppose that $x^*\in \partial f(x)$ and that
$i\in I_\mcF(x)$. 
Then $f_i(x)=f(x)$ and
$(\forall y\in X)$
$f_i(y)-f_i(x)\geq f(y)-f(x)\geq \scal{x^*}{y-x}$.
Therefore, $x^*\in\partial f_i(x)$.
\endproof

\begin{lemma}\label{l:0417a}
Let $(a,b,c)\in X^3$ be colinearly ordered (recall
Definition~\ref{d:colot}). 
Assume \eqref{e:BB} with $I=\{1,2\}$, 
that $\mcF$ is compatible system of functions (recall
Definition~\ref{d:fcomp}), 
that $\DM_{f_1}=\left[a,b\right]$, 
that $\DM_{f_2}=\left[b,c\right]$, and that 
\begin{equation}\label{e:l0417a-6}
\partial f_1(b)\cap\partial f_2(b)\neq\emp.
\end{equation}
Then $f$ is convex and 
\begin{equation}
\label{e:l0417a-4}
(\forall x\in \DM_f)\quad
\partial f(x)=\bigcap_{i\in I_\mcF(x)} \partial f_i(x) = 
\begin{cases}
\partial f_1(x),& \text{if } x\in\left[a,b\right[;\\
\partial f_1(b)\cap\partial f_2(b),&\text{if } x=b;\\
\partial f_2(x),&\text{if } x\in\left]b,c\right].\\
\end{cases}
\end{equation}
\end{lemma}
\proof 
We assume that $a,b,c$ are pairwise distinct since 
the other cases are trivial. First, we show that 
\begin{equation}\label{e:l0417a-2}
\big(\forall x\in\left[a,b\right[\big)
\quad \partial f_1(x)\subseteq\partial f(x),
\end{equation}
Suppose that $x\in\left[a,b\right[$ and that $x^*\in\partial f_1(x)$. 
To establish \eqref{e:l0417a-2}, it suffices to show that 
\begin{equation}\label{e:l0417a-1}
\big(\forall y\in\left[a,c\right]\big)
\quad\scal{x^*}{y-x}\leq f(y)-f(x).
\end{equation}
Indeed, \eqref{e:l0417a-1} is true for $y\in\left[a,b\right]$ by definition of $\partial f_1(x)$ and $f$. 
Now suppose that $y\in\left]b,c\right]$.
By \eqref{e:l0417a-6}, there exists 
$b^*\in\partial f_1(b)\cap\partial f_2(b)$. 
Then 
\begin{subequations}
\begin{align}
\scal{x^*}{y-x}&=\scal{x^*}{b-x}+\scal{x^*}{y-b}\\
&\leq f_1(b)-f_1(x)+
\tfrac{\|y-b\|}{\|b-x\|}\scal{x^*}{b-x}\\
&\leq f_1(b)-f_1(x)+
\tfrac{\|y-b\|}{\|b-x\|}\scal{b^*}{b-x}\\
&\leq f_1(b)-f_1(x)+\scal{b^*}{y-b}\\
&\leq f_1(b)-f_1(x)+f_2(y)-f_2(b)\\
&=f_2(y)-f_1(x)\\
&=f(y)-f(x).
\end{align}
\end{subequations}
Hence \eqref{e:l0417a-1} holds, as does \eqref{e:l0417a-2}.

Switching the roles of $f_1$ and $f_2$, we obtain analogously
\begin{equation}\label{e:l0417a-3}
(\forall x\in\left[c,b\right[\,)\quad \partial f_2(x)\subseteq\partial f(x).
\end{equation}

Next, it is straightforward to check that 
\begin{equation}\label{e:l0417a-5}
\partial f_1(b)\cap\partial f_2(b)\subseteq\partial f(b).
\end{equation}
Since the reverse inclusions of \eqref{e:l0417a-2}, \eqref{e:l0417a-3}, 
and \eqref{e:l0417a-5} follow from Lemma~\ref{l:0502a}, 
we conclude that \eqref{e:l0417a-4} holds.
Using \eqref{e:l0417a-6}, \eqref{e:l0417a-4} and 
Fact~\ref{f:0605a}, we conclude that $\partial f(x)\neq\emp$ 
for all $x\in\left]a,c\right[=\reli \DM_f$.
Finally, it follows from Lemma~\ref{l:0416a} that $f$ is convex.
\endproof

\begin{theorem}\label{t:0513b}
Let $(x_0,\ldots,x_n)\in X^{n+1}$ be colinearly ordered 
(recall Definition~\ref{d:colot}). 
Assume \eqref{e:BB} with $I=\{1,\ldots,n\}$ 
and that $\mcF$ is a compatible system of functions
(recall Definition~\ref{d:fcomp})
such that the following hold:

\begin{enumerate}
\item $(\forall i\in\{1,\ldots,n\})$ $\DM_{f_i}=[x_{i-1},x_i]$.
\item $(\forall i\in\{1,\ldots,n-1\})$
$\partial f_i(x_{i})\cap\partial f_{i+1}(x_i)\neq\emp$.
\end{enumerate}
Then $f$ is convex and 
\begin{equation}
(\forall x\in  \DM_f)\quad \partial f(x)=\bigcap_{i\in I_\mcF(x)}\partial f_i(x);
\end{equation}
\end{theorem}
\proof If $n=1$, then the result is trivial. For $n\geq 2$, 
the result follows by inductively applying Lemma~\ref{l:0417a}.
\endproof

\section{Main results}

\label{s:main}

We are now ready for our main results.

\begin{theorem}[main result I]
\label{t:main1}
Assume \eqref{e:BB}, 
that $\mcF$ is a compatible system of functions 
(recall Definition~\ref{d:fcomp}), 
and that the following hold:
\begin{enumerate}[label={\rm(\alph*)}]
\item\label{t:m1a}
$\DM_f=\bigcup_{i\in I} \DM_{f_i}$ \ is convex and 
at least $2$-dimensional.
\item\label{t:m1b} $\{ \DM_{f_i}\}_{i\in I}$ \ is a compatible
system of sets (recall Definition~\ref{d:comp}).
\item\label{t:m1c} 
There exists a finite subset $E$ of $X$ 
such that 
\begin{equation}\label{e:tm1a}
\left.
\begin{matrix}
x\in (\reli \DM_f)\smallsetminus E\\
\card I(x) \geq 2
\end{matrix}
\right\}
\quad\Rightarrow\quad
\bigcap_{i\in I(x)}\partial f_i(x)\neq\emp.
\end{equation}
\end{enumerate}
Then $f$ is convex and 
\begin{equation}\label{e:tm1b}
(\forall x\in \reli \DM_f)\quad
\emp\neq\partial f(x)\subseteq\bigcap_{i\in I(x)}\partial f_i(x).
\end{equation}
\end{theorem}
\proof
Let $[x,y]\subseteq (\reli \DM_f)\smallsetminus E$. 
By the compatibility in \ref{t:m1b} and Proposition~\ref{p:colin}, 
there exist a colinearly ordered tuple $(x_0,\ldots,x_n)\in
X^{n+1}$ and functions $f_{i_1},\ldots, f_{i_n}$ in $\mcF$ such that
\begin{equation}
x_0=x;\quad x_n=y;\quad\text{and}\quad
\big(\forall k\in\{1,\ldots,n\}\big)\;\;  
[x_{k-1},x_{k}]\subseteq  \DM_{f_{i_k}}\cap\reli \DM_f.
\end{equation}
Define 
\begin{equation}
\big(\forall k\in \{1,\ldots,n\}\big)\quad
g_k:=f_{i_k}+\iota_{[x_{k-1},x_k]}
=f+\iota_{[x_{k-1},x_k]}.
\end{equation}
Using \eqref{e:tm1a}, we see that, for every 
$k\in\{1,\ldots,n-1\}$, 
\begin{equation}
\partial g_k(x_k)\cap\partial g_{k+1}(x_k)\supseteq
\partial f_{i_k}(x_k)\cap\partial f_{i_{k+1}}(x_k)\supseteq
\bigcap_{i\in I(x_k)}\partial f_i(x_k)\neq\emp.
\end{equation}
By applying Theorem~\ref{t:0513b} to the system 
$\{g_k\}_{k\in\{1,\ldots,n\}}$, 
we see that $g=\min_{k\in\{1,\ldots,n\}} g_k=f+\iota_{[x,y]}$
is convex and hence so is $f\big|_{[x,y]}$.
In view of Lemma~\ref{l:0513c}, we obtain the convexity of $f$.
Finally, for all $x\in\reli \DM_f$, 
we have $\partial f(x)\neq\emp$ by Fact~\ref{f:0605a}. 
Therefore, \eqref{e:tm1b} follows from Lemma~\ref{l:0502a}.
\endproof

\begin{remark}
We note an interesting feature in Theorem~\ref{t:main1}. In assumption
\ref{t:m1c}, we require the non-emptiness of the subdifferential
intersection \eqref{e:tm1a} at all relative interior points
except for finitely many points. Since our conclusion says that $f$ is
convex, the subdifferential intersection is nonempty at every point
in $\reli \DM_f$ (see Fact~\ref{f:0605a}). That means, in order to
check the convexity of $f$, we are allowed to ignore verifying
\eqref{e:tm1a} at finitely many points in $\reli\DM_f$. This turn
out to be very convenient since checking \eqref{e:tm1a} at certain
points may not be obvious (see, for example, Example~\ref{ex:1}).
\end{remark}

\begin{remark}[compatibility on the system of domains is essential] 
Theorem~\ref{t:main1} fails if the domain compatibility assumption 
\ref{t:m1b} is omitted: Indeed, 
Suppose that $X=\RR^2$, that $I=\{1,2\}$, 
and that 
\begin{equation}
f_1=\iota_{[0,1]\times[0,1]}
\quad\text{and}\quad
f_2=\iota_{\left[-1,0\right[\times[0,1]}+1.
\end{equation}
Then $\mcF$ is a compatible system of functions. 
Even though $\DM_f=[-1,1]\times[0,1]$ is convex, 
$\{\DM_{f_i}\cap\reli \DM_f\}_{i\in I}$ is \emph{not} a 
compatible system of sets. 
So, Theorem~\ref{t:main1}\ref{t:m1b} is violated. 
Clearly, we can check that $f$ is not convex.
\end{remark}


\begin{theorem}[main result II]
\label{t:main2} 
Assume \eqref{e:BB}, 
that $\mcF$ is a compatible system of functions 
(recall Definition~\ref{d:fcomp}), 
that each $f_i$ is differentiable on  $\inte\DM_{f_i}\neq\emp$, 
and that the following hold:
\begin{enumerate}[label={\rm(\alph*)}]
\item\label{t:m2a0} $\DM_f=\bigcup_{i\in I}\DM_{f_i}$ is convex
and at least $2$-dimensional.
\item\label{t:m2a1} $\{\DM_{f_i}\}_{i\in I}$ is a compatible
system of sets (recall Definition~\ref{d:comp}).
\item\label{t:m2b} 
There exists a finite subset $E$ of $X$ such that
\begin{equation}\label{e:tm2-1}
\left.
\begin{matrix}
x\in (\inte \DM_f)\smallsetminus E\\
\{i,j\}\subseteq I(x)
\end{matrix}
\right\}
\quad\Rightarrow\quad
\lim_{z\to x\atop z\in\inte\DM_{f_i}}\nabla f_i(z)=\lim_{z\to
x\atop z\in\inte\DM_{f_j}}\nabla f_j(z)\quad\text{exists}.
\end{equation}
\end{enumerate}
Then $f$ is convex; moreover, it is continuously differentiable on 
$(\inte \DM_f)\smallsetminus E$.
\end{theorem}
\proof
We will prove the convexity of $f$ by using Theorem~\ref{t:main1}.
Note that it suffices to verify assumption \ref{t:m1c} of 
Theorem~\ref{t:main1}. 
To this end, let $x\in \inte(\DM_f)\smallsetminus
E$ such that $\card I(x)\geq 2$ and denote 
by $u^*_x$ the limit in \eqref{e:tm2-1}. 
Fact~\ref{f:r25-6} and Lemma~\ref{l:0502a} imply
\begin{equation}
u^*_x \in \bigcap_{i\in I(x)}\partial f_i(x).
\end{equation}
So assumption \ref{t:m1c} in Theorem~\ref{t:main1} holds.
Thus, we conclude that $f$ is convex.

Turning towards the differentiability statement, 
let $\O$ be the set of points at which $f$ is differentiable. 
Then $\bigcup_{i\in I}\inte\DM_{f_i} \subseteq \O \subseteq
\inte\DM_{f}$.

Now let $x\in(\inte D_f)\smallsetminus E$. We consider two cases.

{\em Case~1:} $\card I(x)=1$. Then Proposition~\ref{p:sigle}
implies that $x\in\inte D_{f_i}$ for some $i\in I$, which implies
that $x\in\O$. 

{\em Case~2:} $\card I(x)\geq 2$. Since $x\in\inte\DM_f$, we
obtain $N_{\DM_f}(x)=\{0\}$. 
Hence, by Fact~\ref{f:r25-6}, we have
\begin{equation}\label{e:tm2-2}
\partial f(x)=\closu\conv\menge{\lim_{\nnn}
\nabla f(z_n)}{\O\ni z_n\to x}
\end{equation}
Let $(z_n)_\nnn$ be a sequence in $\O$ such that $z_n\to x$,
and let $(\ve_n)_\nnn$ be in $\RPP$ such that $\ve_n\to 0^+$. 
By Fact~\ref{f:r25-5}, $\nabla f\big|_\O$ is the continuous and 
because $\bigcup_{i\in I}\inte\DM_{f_i}$ is dense in $\DM_{f}$, 
there exists a sequence $(y_n)_\nnn$ in 
$\bigcup_{i\in I}\inte D_{f_i}$ such that 
\begin{equation}
(\forall \nnn)\quad
\|y_n-z_n\|\leq\ve_n\quad\text{and}\quad
\|\nabla f(y_n)-\nabla f(z_n)\|\leq\ve_n.
\end{equation}
Combining with \eqref{e:tm2-1}, we deduce that 
\begin{equation}
y_n\to x\quad\text{and}\quad
\lim_{\nnn}\nabla f(z_n)=\lim_{\nnn}\nabla f(y_n)=u^*_x.
\end{equation}
So, \eqref{e:tm2-2} becomes $\partial f(x)=\{u^*_x\}$. 
Thus, $f$ is differentiable at $x$ by Fact~\ref{f:r25-1}.
\endproof

\begin{corollary}
\label{c:main2} 
Assume \eqref{e:BB}, 
that $\mcF$ is a compatible system of functions 
(recall Definition~\ref{d:fcomp}), 
and that the following hold:
\begin{enumerate}[label={\rm(\alph*)}]
\item
$\DM_f=\bigcup_{i\in I}\DM_{f_i}$ is convex
and at least $2$-dimensional.
\item
$\{\DM_{f_i}\}_{i\in I}$ is a compatible
system of sets (recall Definition~\ref{d:comp}).
\item
$f$ is continuously differentiable on $\inte\DM_f$. 
\end{enumerate}
Then $f$ is convex.
\end{corollary}
\proof
This follows from Theorem~\ref{t:main2} with $E=\varnothing$.
\endproof

\begin{remark}[convex interpolation]
\label{r:ci}
When $X=\RR^2$, then Corollary~\ref{c:main2} is known and of 
importance in the convex interpolation of data;
see, e.g.,
\cite[Theorem~1]{GRANDINE-89},
\cite[Proposition~2.2]{LI-99},
\cite[Theorem~3.1]{SCHUMAKER-11}, \cite[Proposition 5.1]{DAHMEN-91}
and the related \cite[Theorem~3.6]{SCHUMAKER-14}
for further details.
\end{remark}

\begin{corollary}
\label{c:main3} 
Assume \eqref{e:BB}, 
that $\mcF$ is a compatible system of functions 
(recall Definition~\ref{d:fcomp}), 
that each $\DM_{f_i}$ is closed, 
and that 
$(\forall i\in I)(\forall x\in \DM_{f_i})$ 
$f_i(x) = \tfrac{1}{2}\scal{x}{A_ix} + \scal{b_i}{x} + \gamma_i$,
where $A_i\colon X\to X$ is linear with $A_i = A_i^* \succeq 0$, $b_i\in X$, and $\gamma_i\in\RR$.
Furthermore, assume that 
$\DM_f=\bigcup_{i\in I}\DM_{f_i}$ is convex
and at least $2$-dimensional, 
and that 
\begin{equation}
\left.
\begin{matrix}
\{i,j\} \subseteq I\\
i\neq j\\
x\in \DM_{f_i}\cap \DM_{f_j}
\end{matrix}
\right\}
\quad\Rightarrow\quad
A_ix+b_i = A_jx+b_j.
\end{equation}
Then $f$ is convex; moreover, it is continuously differentiable on 
$\inte\DM_{f}$. 

\end{corollary}
\proof
This follows from Example~\ref{ex:0815} and Theorem~\ref{t:main2} with $E=\varnothing$ since
$\nabla f_i(x) = A_ix+b_i$ when $x\in\inte\DM_{f_i}$. 
\endproof

\begin{remark}[piecewise linear-quadratic function]
\label{r:plq}
Consider Corollary~\ref{c:main3} with the additional assumption
that each $\DM_{f_i}$ is a polyhedral set.
Then Corollary~\ref{c:main3} provides a sufficient condition
for checking the convexity of $f$ which in this case is a
\emph{piecewise linear-quadratic function}.
These functions play a role in computer-aided convex analysis
(see \cite{LUCET-10} and also 
\cite[Section~10.E]{ROCKAFELLAR-98a} for further information
on functions of this type). 
Moreover, we thus partially answer an open question from
\cite[Section 23.4.2]{LUCET-12}. Our results also enhance our
understanding of how nonconvexity occurs and form another step
towards building a nonconvex toolbox that extends current bivariate
computational convex analysis
algorithms~\cite{GARDINER-13,GARDINER-11}. 
\end{remark}

\section{Checking convexity}

\label{s:Checkconv}

We start with an application of Theorem~\ref{t:main2}.

\begin{example}\label{ex:1}
The function 
\begin{equation}\label{e:0319a}
f\colon \RR^2\to\RR \colon (x,y) \mapsto
\begin{cases}
\disp\frac{x^2+y^2+2\max\{0,xy\}}{|x|+|y|},&\text{if}\quad(x,y)\neq (0,0);\\
0,&\text{if}\quad(x,y)=(0,0)
\end{cases}
\end{equation}
is convex, and differentiable on $\RR^2\smallsetminus\{(0,0)\}$. 
\end{example}
\proof
First, set $I:=\{1,\ldots,4\}$ and
\begin{subequations}\label{e:0320a}
\begin{align}
f_1(x,y)&:=
	\begin{cases}
	x+y,&\text{if}\quad (x,y)\in A_1:=\RP^2,\\
	+\infty,&\text{otherwise};
	\end{cases}\\
f_2(x,y)&:=
	\begin{cases}
	\frac{x^2+y^2}{-x+y},&\text{if}\quad (x,y)\in A_2:=\RM\times\RP,\\
	+\infty,&\text{otherwise};
	\end{cases}\\
f_3(x,y)&:=
	\begin{cases}
	-x-y,&\text{if}\quad (x,y)\in A_3:=\RM^2,\\
	+\infty,&\text{otherwise};
	\end{cases}\\
f_4(x,y)&:=
	\begin{cases}
	\frac{x^2+y^2}{x-y},&\text{if}\quad (x,y)\in A_4:=\RP\times\RM,\\
	+\infty,&\text{otherwise}.
	\end{cases}
\end{align}
\end{subequations}
Then $\{f_i\}_{i\in I}$ is a compatible system of functions (recall 
Definition~\ref{d:fcomp}) with $f$ being the corresponding
piecewise-defined
function. Moreover, $\{D_{f_i}\}_{i\in I}=\{A_i\}_{i\in I}$ is a
compatible system of sets.

Each $f_i$ is differentiable on $\inte A_i$ with the gradient given by
\begin{subequations}\label{e:0319b}
\begin{align}
\nabla f_1(x,y)&=(1,1)\quad\text{for}\quad (x,y)\in\inte A_1;\\
\nabla f_2(x,y)&=
\Big(\tfrac{-x^2+2xy+y^2}{(x-y)^2},\tfrac{-x^2-2xy+y^2}{(x-y)^2}\Big)
\quad\text{for}\quad (x,y)\in\inte A_2;\\
\nabla f_3(x,y)&=(-1,-1)
\quad\text{for}\quad (x,y)\in\inte A_3;\\
\nabla f_4(x,y)&=
\Big(\tfrac{x^2-2xy-y^2}{(x-y)^2},\tfrac{x^2+2xy-y^2}{(x-y)^2}\Big)
\quad\text{for}\quad (x,y)\in\inte A_4.
\end{align}
\end{subequations}
One readily checks that the Hessian of each $f_i$ is positive
semi-definite on $\inte A_i$;
hence, by the continuity of $f_i\big|_{\DM_{f_i}}$, 
we have that $f_i$ is convex.

Moreover, 
\begin{subequations}\label{e:0320g}
\begin{align}
&(\forall a>0)\quad
\lim_{(x,y)\to(a,0)\atop (x,y)\in\inte A_1} \nabla f_1(x,y)
=\lim_{(x,y)\to(a,0)\atop (x,y)\in\inte A_4} \nabla f_4(x,y)=(1,1);\\
&(\forall a<0)\quad
\lim_{(x,y)\to(a,0)\atop (x,y)\in\inte A_2} \nabla f_2(x,y)
=\lim_{(x,y)\to(a,0)\atop (x,y)\in\inte A_3} \nabla f_3(x,y)=(-1,-1);\\
&(\forall b>0)\quad
\lim_{(x,y)\to(0,b)\atop (x,y)\in\inte A_1} \nabla f_1(x,y)
=\lim_{(x,y)\to(0,b)\atop (x,y)\in\inte A_2} \nabla f_2(x,y)=(1,1);\\
&(\forall b<0)\quad
\lim_{(x,y)\to(0,b)\atop (x,y)\in\inte A_3} \nabla f_3(x,y)
=\lim_{(x,y)\to(0,b)\atop (x,y)\in\inte A_4} \nabla f_4(x,y)=(-1,-1).
\end{align}
\end{subequations}
Now set $E:=\{(0,0)\}$. From the above
computations, we observe that all assumptions of Theorem~\ref{t:main2} are satisfied.
Thus, we conclude that $f$ is a convex function that is
also continuously differentiable away from the origin. 
\endproof

In fact, the function defined by \eqref{e:0319a} is actually a norm since
it is clearly positively homogeneous. 
An analogous use of Theorem~\ref{t:main2} allows for a systematic proof of the convexity of the
function considered next. 

\begin{example}\label{ex:4}
The function 
\begin{equation}
f\colon \RR^2\to\RR\colon (x,y)\mapsto \begin{cases}
\sqrt{x^6+y^4},&\text{if}\quad xy\geq 0;\\
|x|^3+y^2,&\text{otherwise},
\end{cases}
\end{equation}
is a convex and continuously differentiable. 
\end{example}

We conclude this section with an application of Theorem~\ref{t:main1}.

\begin{example}\label{ex:5}
The function 
\begin{equation}
f\colon \RR^2\to\RR\colon (x,y)\mapsto 
f(x,y):=\begin{cases}
\sqrt{x^4+y^2},&\text{if}\quad xy \geq 0;\\
x^2+|y|,&\text{otherwise},
\end{cases}
\end{equation}
is convex. 
\end{example}
\proof
First, set $I:=\{1,\ldots,4\}$ and
\begin{subequations}
\begin{align}
f_1(x)&:=
	\begin{cases}
	\sqrt{x_1^4+x_2^2},&\text{if}\quad x\in A_1:=\RP^2,\\
	+\infty,&\text{otherwise};
	\end{cases}\\
f_2(x)&:=
	\begin{cases}
	x_1^2+x_2,&\text{if}\quad x\in A_2:=\RM\times\RP,\\
	+\infty,&\text{otherwise};
	\end{cases}\\
f_3(x)&:=
	\begin{cases}
	\sqrt{x_1^4+x_2^2},&\text{if}\quad x\in A_3:=\RM^2,\\
	+\infty,&\text{otherwise};
	\end{cases}\\
f_4(x)&:=
	\begin{cases}
	x_1^2-x_2,&\text{if}\quad x\in A_4:=\RP\times\RM,\\
	+\infty,&\text{otherwise}.
	\end{cases}
\end{align}
\end{subequations}
Then $\{f_i\}_{i\in I}$ is a compatible system of functions (recall 
Definition~\ref{d:fcomp}) with $f$ being the corresponding piecewise
function. Moreover, $\{D_{f_i}\}_{i\in I}=\{A_i\}_{i\in I}$ is a
compatible system of sets.

Each $f_i$ is differentiable on $\inte A_i$ with the gradient given by
\begin{subequations}\label{e:ex5a}
\begin{align}
\nabla f_1(x)&=\Big(\tfrac{2x_1^3}{\sqrt{x_1^4+x_2^2}},
\tfrac{x_2}{\sqrt{x_1^4+x_2^2}}\Big)
\quad\text{for}\quad x\in\inte A_1;\\
\nabla f_2(x)&=(2x_1,1)\quad\text{for}\quad x\in\inte A_2;\\
\nabla f_3(x)&=\Big(\tfrac{2x_1^3}{\sqrt{x_1^4+x_2^2}},
\tfrac{x_2}{\sqrt{x_1^4+x_2^2}}\Big)
\quad\text{for}\quad x\in\inte A_3;\\
\nabla f_4(x)&=(2x_1,-1)\quad\text{for}\quad x\in\inte A_4.
\end{align}
\end{subequations}

Next, since the Hessian of $f_i$ is positive semidefinite on $\inte A_i$, 
we deduce that each $f_i$ is convex.

Now set $E:=\{(0,0)\}$. We will verify 
\eqref{e:tm1a}. Note that simple computations show the following:

For $x=(0,x_2)\in (A_1\cap A_2)\smallsetminus E$,
\begin{equation}
\lim_{z\to x\atop z\in\inte A_1}\nabla f_1(z)=
\lim_{z\to x\atop z\in\inte A_2}\nabla f_2(z)=(0,1)
\in\partial f_1(x)\cap\partial f_2(x).
\end{equation}

For $x=(0,x_2)\in (A_3\cap A_4)\smallsetminus E$,
\begin{equation}
\lim_{z\to x\atop z\in\inte A_3}\nabla f_3(z)=
\lim_{z\to x\atop z\in\inte A_4}\nabla f_4(z)=(0,-1)
\in\partial f_3(x)\cap\partial f_4(x).
\end{equation}

For $x=(x_1,0)\in (A_2\cap A_3)\smallsetminus E$,
\begin{subequations}
\begin{align}
&\lim_{z\to x\atop z\in\inte A_2}\nabla f_2(z)=(2x_1,1)
\quad\text{and}\quad
N_{A_2}(x)=\{0\}\times\RM;\\
&\lim_{z\to x\atop z\in\inte A_3}\nabla f_3(z)=(2x_1,0)
\quad\text{and}\quad
N_{A_3}(x)=\{0\}\times \RP. 
\end{align}
\end{subequations}
Then, using Fact~\ref{f:r25-6}, we conclude that
$\partial f_2(x)\cap\partial f_3(x)\neq\emp$.

For $x=(x_1,0)\in (A_1\cap A_4)\smallsetminus E$,
\begin{subequations}
\begin{align}
&\lim_{z\to x\atop z\in\inte A_1}\nabla f_1(z)=(2x_1,0)
\quad\text{and}\quad
N_{A_1}(x)=\{0\}\times\RM;\\
&\lim_{z\to x\atop z\in\inte A_4}\nabla f_4(z)=(2x_1,-1)
\quad\text{and}\quad
N_{A_4}(x)=\{0\}\times\RP.
\end{align}
\end{subequations}
Then, using Fact~\ref{f:r25-6}, we conclude that
$\partial f_1(x)\cap\partial f_4(x)\neq\emp$.

So, we have verified that assumption \ref{t:m1c} in Theorem~\ref{t:main1} holds.
Therefore, $f$ is convex by Theorem~\ref{t:main1}.
\endproof

\section{Detecting the lack of convexity}

\label{s:Lackconv}

The nonempty subdifferential intersection condition \eqref{e:tm1a}
is indeed crucial for the check of convexity: we will
see in the following result that the violation of \eqref{e:tm1a} leads
to nonconvexity.

\begin{theorem}[detecting lack of convexity]
\label{t:noncvx}
Assume \eqref{e:BB}, that $\mcF$ is 
a compatible system of functions (recall Definition~\ref{d:fcomp}), 
and that 
\begin{equation}
(\exists x\in\reli\DM_f)\quad\bigcap_{i\in I(x)}\partial f_i(x)=\emp.
\end{equation}
Then $f$ is not convex. 
\end{theorem}
\proof Using Lemma~\ref{l:0502a}, we have 
$\partial f(x)\subseteq\bigcap_{i\in I(x)}\partial f_i(x)=\emp$. 
Therefore, by Fact~\ref{f:0605a}, $f$ is not convex.
\endproof

Using Theorem~\ref{t:noncvx}, we will now illustrate that 
the finiteness assumption on $E$ is important for our main results 
(Theorem~\ref{t:main1} and Theorem~\ref{t:main2}).

\begin{example}
\label{ex:3}
Suppose that $X=\RR^2$,
set 
\begin{subequations}
\begin{align}
f_1(x,y)&:=\begin{cases}
\max\{-x,y\},&\text{if}\quad (x,y)\in A_1 := \RP\times\RR;\\
+\infty,&\text{otherwise},
\end{cases}\\
f_2(x,y)&:=\begin{cases}
\max\{x,y\},&\text{if}\quad (x,y)\in A_2 := \RM\times\RR;\\
+\infty,&\text{otherwise},
\end{cases}
\end{align}
\end{subequations}
$\mcF:=\{f_1,f_2\}$, $f := \min\{f_1,f_2\}$, i.e., 
\begin{equation}
f\colon \RR^2\to\RR\colon (x,y)\mapsto \max\big\{ -|x|,y\big\},
\end{equation}
and 
$E:=\{0\}\times\RM$. 
Then one checks the following:
\begin{enumerate}
\item $\mcF$ is a compatible system of functions.
\item $\{A_1,A_2\}$ is a compatible system of sets. 
\item For every $(\ox,\oy)\in \RR^2 \smallsetminus E$ with 
$I_\mcF(\ox,\oy)=\{1,2\}$, we must have
$(\ox,\oy)\in\{0\}\times\RPP$, i.e., $\ox=0$ and $\oy>0$. 
Then $f(x,y)=y$ locally around $(0,\oy)$ and thus
\begin{equation}
\partial f_1(0,\oy)\cap\partial f_2(0,\oy)\supseteq\partial f(0,\oy)=\{(0,1)\}.
\end{equation}
\end{enumerate}
So, all assumptions in Theorems~\ref{t:main1} and \ref{t:main2}
are satisfied except that $E$ is infinite. 
However, for $(0,\oy)\in E\smallsetminus\{(0,0)\}$, we have $\oy<0$; 
thus, $f_1(x,y)=-x+\iota_{A_1}(x,y)$ locally around $(0,\oy)$. 
It follows that
\begin{equation}
\partial f_1(0,\oy)=(-1,0)+N_{A_1}(0,\oy)=
\left]-\infty,-1\right]\times\{0\}
\end{equation}
and similarly that 
\begin{equation}
\partial f_2(0,\oy)=(1,0)+N_{A_2}(x,y)=
\left[1,+\infty\right[\times\{0\}.
\end{equation}
Hence $\partial f_1(0,y)\cap\partial f_2(0,y)=\emp$
and so 
$f$ is not convex by applying Theorem~\ref{t:noncvx} or
by direct inspection. 
\end{example}


In the previous example, the set $E$ was infinite, but unbounded. 
In the next (slightly more involved) example, we provide a case where
$E$ is bounded. 

\begin{example}\label{ex:2}
Suppose that $X=\RR^2$, 
set $I:=\{1,\ldots,6\}$, 
\begin{subequations}
\begin{align}
&A_1:=\menge{(x,y)\in  \RR^2}{|x|+|y|\geq 1,x\geq0,y\geq0},\\
&A_2:=\menge{(x,y)\in  \RR^2}{|x|+|y|\geq 1,x\leq0,y\geq0},\\
&A_3:=\menge{(x,y)\in  \RR^2}{|x|+|y|\geq 1,x\leq0,y\leq0},\\
&A_4:=\menge{(x,y)\in  \RR^2}{|x|+|y|\geq 1,x\geq0,y\leq0},\\
&A_5:=\menge{(x,y)\in  \RR^2}{|x|+|y|\leq 1,x\geq0},\\
&A_6:=\menge{(x,y)\in  \RR^2}{|x|+|y|\leq 1,x\leq0},
\end{align}
\end{subequations}
\begin{equation}
f\colon \RR^2\to \RR\colon 
(x,y)\mapsto 
\max\big\{1-|x|,|y|\big\},
\end{equation}
\begin{equation}
\mcF:=\{f_i\}_{i\in I},\quad\text{where $~(\forall i\in I)\quad 
f_i:=f+\iota_{A_i}$,}
\end{equation}
and 
$E:=\big(\{0\}\times\left[-1,1\right]\big)\cup\{(\pm 1,0)\}$. 
Then one checks the following:
\begin{enumerate}
\item Each $f_i$ is convex and continuous on $\DM_{f_i}$ because
\begin{subequations}
\begin{align}
f_i(x,y)&=y+\iota_{A_i}(x,y)\quad\text{for}\quad i\in\{1,2\};\\
f_i(x,y)&=-y+\iota_{A_i}(x,y)\quad\text{for}\quad i\in\{3,4\};\\
f_5(x,y)&=-x+1+\iota_{A_5}(x,y);\quad\text{and}\\
f_6(x,y)&=x+1+\iota_{A_6}(x,y).
\end{align}
\end{subequations}
Consequently, $\mcF$ is a compatible system of functions. 
\item $\{A_i\}_{i\in I}$ is a compatible system of sets.
\item $f$ is the piecewise-defined function associated with
$\mcF$.
\item Take $(\ox,\oy)\in \RR^2\smallsetminus E$ with $\card I_\mcF(\ox,\oy)\geq 2$. 
Then
\begin{equation}
I_\mcF(\ox,\oy)\in
\big\{\{1,2\},\{1,4\},\{1,5\},\{2,3\},\{2,6\},\{3,4\},\{3,6\},
\{4,5\}\big\}.
\end{equation}
Suppose, for instance, that 
$I_\mcF(\ox,\oy)=\{1,5\}$. Then $\ox>0$, $\oy>0$, and $\ox+\oy=1$. We have
\begin{equation}
\partial f_1(\ox,\oy)=(0,1)+N_{A_1}(\ox,\oy)=(0,1)+\RP(-1,-1)
\end{equation}
and
\begin{equation}
\partial f_5(\ox,\oy)=(0,1)+N_{A_5}(\ox,\oy)=(-1,0)+\RP(1,1).
\end{equation}
Then 
\begin{equation}
(-\tfrac{1}{2},\tfrac{1}{2})\in\partial f_1(\ox,\oy)\cap\partial f_5(\ox,\oy);
\end{equation}
similarly, one obtains nonemptiness for the other cases. 
\end{enumerate}
We observe that all assumptions in Theorem~\ref{t:main1} and 
Theorem~\ref{t:main2} are satisfied except that $E$ is infinite and bounded.
However, for every $(0,\oy)\in\{0\}\times\left]-1,1\right[\subseteq E$, 
we have $f(x,y)=\min\{f_5(x,y),f_6(x,y)\}=-|x|+1$ 
locally around $(0,\oy)$. Clearly, $f$ is not convex.
\end{example}


\section*{Acknowledgments}
HHB was partially supported by the Natural Sciences and Engineering
Research Council (NSERC) of Canada and by the Canada Research Chair Program.
YL was partially supported by the Natural Sciences and Engineering Research
Council of Canada through a Discovery grant.
HMP was partially supported by an NSERC accelerator grant of HHB.



\begin{thebibliography}{999}

\sepp

\bibitem{BC2011}
H.H.~Bauschke and P.L.~Combettes,
\emph{Convex Analysis and Monotone Operator Theory in Hilbert Spaces},
Springer, 2011.

\bibitem{DAHMEN-91}
W.~Dahmen,
\emph{Convexity and {B}ernstein-{B}\'ezier polynomials}, in
{\em Curves and surfaces ({C}hamonix-{M}ont-{B}lanc, 1990)}, 
pages~107--134, Academic Press, 1991. 

\bibitem{GARDINER-13}
B.~Gardiner, K.~Jakee, and Y.~Lucet,
Computing the partial conjugate of convex piecewise linear-quadratic
  bivariate functions,
\emph{Computational Optimization and Applications}~58 (2014),
249--272. 

\bibitem{GARDINER-11}
B.~Gardiner and Y.~Lucet,
Computing the conjugate of convex piecewise linear-quadratic
  bivariate functions,
{\em Mathematical Programming (Series B)}~139 (2013), 161--184.

\bibitem{GRANDINE-89}
T.A.~Grandine,
On convexity of piecewise polynomial functions on triangulations,
\emph{Computer Aided Geometric Design}~6 (1989), 181--187.

\bibitem{LI-99}
A.~Li,
Convexity preserving interpolation,
\emph{Computer Aided Geometric Design}~16 (1999), 127--147.

\bibitem{LUCET-10}
Y.~Lucet,
What shape is your conjugate? {A} survey of computational convex
  analysis and its applications,
\emph{SIAM Review}~52 (2010), 505--542.

\bibitem{LUCET-12}
Y.~Lucet,
Techniques and open questions in computational convex analysis,
in {\em Computational and Analytical Mathematics}, pages 485--500, 
Springer, 2013.

\bibitem{MorNam}
B.S.~Mordukhovich and N.M.~Nam,
\emph{An easy path to convex analysis and applications},
 Morgan \& Claypool, 2014. 

\bibitem{ROCKAFELLAR-70a}
R.T.~Rockafellar,
{\em Convex Analysis},
Princeton University Press, 1970.

\bibitem{ROCKAFELLAR-98a}
R.T.~Rockafellar and R.J-B.~Wets,
{\em Variational Analysis},
Springer, 1998.

\bibitem{SCHUMAKER-11}
L.L. Schumaker and H.~Speleers,
Convexity preserving splines over triangulations,
\emph{Computer Aided Geometric Design}~28 (2011), 270--284. 

\bibitem{SCHUMAKER-14}
L.L.~Schumaker and H.~Speleers,
Convexity preserving {$C^0$} splines,
\emph{Advances in  Computational Mathematics}~40 (2014), 117--135.

\bibitem{ZALINESCU-02}
C.~Z{\u{a}}linescu,
{\em Convex Analysis in General Vector Spaces},
World Scientific, 2002.

\end{thebibliography}

\end{document}